\documentclass[12pt]{article}
\usepackage{amsmath}
\usepackage{color}
\usepackage{fancyhdr}
\usepackage{verbatim}
\usepackage{amsfonts,longtable,mathrsfs}
\usepackage{epsfig}
\usepackage{amsfonts}
\usepackage{color}
\usepackage{fullpage}
\usepackage[numbers,sort&compress]{natbib}

\def\ec{\end{center}}
\def\bc{\begin{center}}
\def\ec{\end{center}}
\newtheorem{definition}{Definition}[section]
\newtheorem{theorem}{Theorem}[section]

\begin{document}
\bc {\bf On certain rational recursive sequences of order four
  }\ec
\medskip
\bc
Mensah Folly-Gbetoula\footnote{Corresponding author:\\ Mensah.Folly-Gbetoula@wits.ac.za(M. Folly-Gbetoula)\\ 
ORCID: 0000-0002-3046-0679}
and Darlison Nyirenda\footnote{Author: \\Darlison.Nyirenda@wits.ac.za(D. Nyirenda)\\
}
\vspace{1cm}
\\School of Mathematics, University of the Witwatersrand, Wits 2050, Johannesburg, South Africa.\\

\ec
\begin{abstract}
\noindent We obtain solutions to the recursive sequences of the form
$$x_{n + 1} = \frac{x_{n - 3}x_{n }}{x_{n - 2}(a_n + b_nx_{n -3}x_{n})}$$
where $a_n$ and $b_n$ are arbitrary sequences of real numbers, and the initial values are gives as; $x_{-3},x_{-2}, x_{-1}$ and $x_{0}$. Our methodology is to employ a group theoretic method which lowers the order of the equations, and then solve the resulting lower order recurrence relations that arise therefrom.
\end{abstract}
\section{Introduction} \setcounter{equation}{0}
Due to the work of Sophus Lie on differential equations \cite{Lie}, there has been great interest in the invariance properties of differential equations under groups of point transformations. Lie discovered that most of  the known theories for solving differential equations are closely related to the idea of infinitesimal transformations. One of the applications of symmetry analysis (Lie's method) in differential equations is the obtention of solutions. A lot of work has been published that reflects the application of the method to difference equations \cite{FK2, FK,hydon0, hydon1, JV, W, FN}. Symmetry analysis on difference equations started with Shigeru Maeda \cite{Maeda, Maeda2} who studied the continuous point symmetries of difference equations. Nalini Joshi \cite{JV} came up with a concrete method for finding symmetries of first order difference equations and gave the local analytic diffeomorphism that makes first order discrete dynamical systems linear. In the recent past, Hydon \cite{hydon0, hydon1} came up with algorithms for finding symmetries and first integrals of difference equations independent of their order. It must be stated that, applying these methods involves a lot of cumbersome and tedious calculations. \par \noindent
In this paper, we perform a symmetry analysis and derive closed form formulas for solutions of difference equations of the form
\begin{align}\label{xn}
x_{n+1}=\frac{x_{n - 3}x_{n }}{x_{n - 2}(a_n + b_nx_{n -3}x_{n})},
\end{align}
where $a_n$ and $b_n$ are sequences of real numbers. For related work on this, see \cite{bc,ibm0,FK2, FN, EAA,KE}.
\subsection{Preliminaries}
We start by recalling some terminology in Lie analysis of differential and difference equations mostly taken from \cite{BA,hydon0}.
\begin{definition}\cite{BA}
Let $x = (x_1; x_2; \dots ; x_n)$ lie in a region $D\subset \mathcal{R}^n$. The
set of transformations
\begin{align}\label{transfo}
\tilde{x} = X(x; \varepsilon),
\end{align}
defined for each $x$ in $D$ and parameter $\varepsilon$ in a set $\mathcal{S}$, a subset of $\mathcal{R}$, with $\phi(\varepsilon, \delta)$ defining a law of composition of parameters $\varepsilon$ and $\delta$ in $\mathcal{S}$, constitutes a one-parameter group of transformations on $D$ if the following hold:
\begin{itemize}
\item
For each $\varepsilon$ in $\mathcal{S}$, the transformation is one-to-one  and onto. (Hence, $\tilde{x}$  lies in $D$).
\item
 $S$ with the law of composition $\phi$ forms a group $G$.
\item
For each $x$ in $D$, $\tilde{x} = x$, when $\varepsilon=\varepsilon _0$ corresponds to the identity $e$, that is,
\begin{align}
X(x; \varepsilon _0) = x.
\end{align}
\item
If $\tilde{x} = X(x;\varepsilon), \tilde{\tilde{x}} = X(\tilde{x}; \delta)$, it follows that
\begin{align}
\tilde{\tilde{x}} =X(x;\phi(\varepsilon, \delta)).
\end{align}
\end{itemize}
\end{definition}
The above definition together with the following:
\begin{itemize}
\item
$\varepsilon$ is a continuous parameter, that is, $\mathcal{S}$ is an interval in $\mathcal{R}$. Without a loss of generality, $\varepsilon = 0$ corresponds to the identity element $e$:
\item
 $X$ is infinitely differentiable with respect to $x$ in $D$ and is an analytic function of $\varepsilon$ in $\mathcal{S}$.
 \item
$\phi(\varepsilon, \delta)$ is an analytic function of $\varepsilon$ and $\delta$; $\varepsilon \in \mathcal{S}$; $\delta \in \mathcal{S}$,
\end{itemize}
define a one-parameter group of transformations.
\begin{definition}\cite{BA}
The infinitesimal generator of the one-parameter Lie group of transformations \eqref{transfo} is the operator
\begin{align}
X=X(x)=\xi(x)\times \Delta = \sum_{s=1}^n\xi _s (x)\frac{\partial}{\partial x_s},
\end{align}
where $\Delta$ is the gradient operator
\begin{align}
\Delta = \left( \frac{\partial}{\partial x_1},\frac{\partial}{\partial x_2},\dots,\frac{\partial}{\partial x_n}\right).
\end{align}
\end{definition}
\begin{definition}\cite{BA}
An infinitely differentiable function $F(x)$ is an invariant function of the Lie group of transformations \eqref{transfo} if and only if, for any group transformations,
\begin{align}
F(\tilde{x}) = F(x).
\end{align}
\end{definition}
\begin{theorem}\cite{BA}
$F(x)$ is invariant under the Lie group of transformations \eqref{transfo} if and only if,
\begin{align}
XF(x) = 0.
\end{align}
\end{theorem}
Considering the $r$th-order difference equation
\begin{align}\label{1.1}
u_{n+ r+1}=\varphi (n, u_n, u_{n+1},\dots, u_{n+ r}).
\end{align}
We employ the infinitesimal group of point transformations
\begin{align}
(n, u_{n+s}) \mapsto (\tilde{n}=n, \tilde{u}_{n+s} = u_{n+s} + \varepsilon \xi (n+s,u_{n+s}))
\end{align}
and let
\begin{align}
\mathcal{X}=\xi (n,u_n)\frac{\partial}{\partial u_n}
\end{align}
be the generator of the group transformations with prolonged form
\begin{align}
X= \mathcal{X}^{[r]}=\xi (n,u_n)\frac{\partial\quad}{\partial u_n}+\xi (n+1,u_{n+1})\frac{\partial\quad}{\partial u_n}+\dots+\xi (n+r,u_{n+r})\frac{\partial\qquad}{\partial u_{n+r}}.
\end{align}
A symmetry is a group of transformations that maps solutions onto other solutions. Here, the criterion of invariance is given by
\begin{align}\label{1.1transfo}
\tilde{u}_{n+r+1}=\varphi (\tilde{n}, \tilde{u}_n, \tilde{u}_{n+1},\dots, \tilde{u}_{n+r}),
\end{align}
whenever \eqref{1.1} holds. As a result, the linearized symmetry condition becomes
\begin{align}\label{LSC}
S^ {r+1}\xi(n,u_{n})-X\varphi =0.
\end{align}
\section{Main results}
We consider ordinary difference equations;
\begin{align}\label{un}
u_{n+4}=\varphi=\frac{u_{n}u_{n+3}}{u_{n+1}(A_n+ B_n u_{n}u_{n+3})},
\end{align}
where $A_n$ and $B_n$ are random real sequences, equivalent to \eqref{xn}. \par \noindent
By the criterion of invariance \eqref{LSC},
\begin{align}\label{a1}
 &\xi(n+4,\varphi)-\frac{A_nu_n\xi(n+3,u_{n+3})}{u_{n+1}(A_n +B_nu_nu_{n+3})^2}+\frac{u_nu_{n+3}\xi(n+1,u_{n+1})}{u_{n+1}^2(A_n +B_nu_nu_{n+3})}\nonumber\\
 &-\frac{A_nu_{n+3}\xi(n,u_{n})}{u_{n+1}(A_n +B_nu_nu_{n+3})^2}=0.
\end{align}
To solve for $\xi$, we differentiate  implicitly, the functional equation \eqref{a1} with respect to $u_n$ by viewing $u_{n+1}$ as a function of $u_n$, $u_{n+3}$ and $\varphi$. That is, to operate
$$L=\frac{\partial\quad }{\partial {u_n}}+\frac{\partial u_{n+1}}{\partial u_{n}\quad}\frac{\partial\qquad }{\partial u_{n+1}}=\frac{\partial\quad }{\partial {u_n}}-\frac{\frac{\partial \varphi}{\partial u_{n}}}{\frac{\partial \varphi\quad}{\partial u_{n+1}}}\frac{\partial\qquad }{\partial u_{n+1}}$$
on \eqref{a1}. This leads to
\begin{align}\label{a3}
& B_n{u_n}^2u_{n+1}\xi(n+3,u_{n+3})+\left(A_n{u_n} +B_n{u_n}^2u_{n+3}\right)u_{n+1}\xi'(n+1,u_{n+1})
-(A_nu_n \nonumber\\
&+B_n{u_n}^2u_{n+3}){\xi(n+1,u_{n+1})}-(A_nu_n +B_n{u_n}^2u_{n+3})u_{n+1}\xi'(n,u_n)+(2B_nu_nu_{n+3}\nonumber\\
&+{A_n})u_{n+1}\xi(n,u_n)=0
\end{align}
after clearing fractions.
By thrice differentiating \eqref{a3} with respect to $u_n$, keeping $u_{n+1}$ fixed, we obtain
\begin{align}\label{a4}
-u_nu_{n+1}\left(A_n +B_nu_nu_{n+3}\right)\xi ^{(4)}(n,u_n)-2u_{n+1}(A_n+2B_nu_nu_{n+3})\xi^{(3)}(n,u_n)=0.
\end{align}
We separate \eqref{a4} based on the fact that the function $\xi$ depends on the continuous variable $u_n$ only. We get
\begin{align}\label{a5}
u_{n+1}u_{n+3} \quad \text{terms} &: u_n \xi^{(4)}(n,u_n)+4\xi ^{(3)}(n,u_n)=0\\
u_{n+1} \qquad \quad\text{terms} &:u_n \xi^{(4)}(n,u_n)+2\xi ^{(3)}(n,u_n)=0.
\end{align}
Equivalently,
\begin{align}\label{a6}
&\xi(n, u_n) = \beta _n u_n ^2 + \gamma _n u_n+\lambda _n
\end{align}
for some functions $\beta _n,\,\gamma _n$ and $\lambda _n$ of $n$.\par \noindent
Here, we perform a substitution of \eqref{a6} into \eqref{a1} and separate the resulting equation by setting the coefficients of products of shifts of $u_n$ to zero (since the functions $\beta _n,\,\gamma _n$ and $\lambda _n$ depend on $n$ only). This leads to the following system:
\begin{subequations}\label{a7}
\begin{align}
{u_n}^2{u_{n+1}}^2{u_{n+3}}^2\; \text{terms}&:\qquad B_n ^2 \lambda _{n+4}+ B_n \beta _{n+1}=0\\
{u_n}{u_{n+1}}^2{u_{n+3}}\; \text{terms}&: \qquad 2A_nB_n\lambda _{n+1}+A_n\beta _{n+1}=0\\
{u_n}^2{u_{n+1}}{u_{n+3}}^2\; \text{terms}&:\qquad B_n\gamma _{n+1}+B_n\gamma _{n+4}=0\\
{u_n}^2{u_{n+1}}{u_{n+3}}\; \text{terms}&:\qquad -A_n\beta _{n}=0\\
{u_n}^2{u_{n+3}}^2\; \text{terms}&:\qquad B_n\lambda _{n+1}+\beta _{n+4}=0\\
{u_n}{u_{n+1}}{u_{n+3}}^2\; \text{terms}&:\qquad -A_n\beta _{n+3}=0\\
{u_{n+1}}^2\; \text{terms}&:\qquad{A_n}^2\lambda _{n+4} =0\\
{u_n}{u_{n+1}}{u_{n+3}}\; \text{terms}&:\qquad -A_n \gamma _n+A_n \gamma _{n+1}-A_n \gamma _{n+3}+A_n \gamma _{n+4}=0\\
{u_n}{u_{n+1}}\; \text{terms}&:\qquad-A_n \lambda _{n+3}=0\\
{u_n}{u_{n+3}}\; \text{terms}&:\qquad A_n \lambda _{n+1}=0\\
{u_{n+1}}{u_{n+3}}\; \text{terms}&: \qquad-A_n \lambda _n =0.
\end{align}
\end{subequations}
Solving \eqref{a7} leads to the following \lq final constraint\rq:
\begin{align}
&\beta _n=0\\& \lambda _n =0\\
&\gamma _{n}+\gamma _{n+3}=0.\label{rel}
\end{align}
It is not difficult to see that $(-1)^n,\; \exp(i\pi n/3)$ and $\exp(-i\pi n/3)$
are the solutions of \eqref{rel}.\par \noindent
Set $\gamma  =\exp\left( i{\pi}/{3}\right)$. Thanks to \eqref{a6}, the three characteristics are thus given as
\begin{subequations}\label{chara}
\begin{align}\label{5charac}
\xi_1=&(-1) ^n u_n ,\; \xi_2=\gamma ^n u_n \;\text{ and }\;
 {\xi}_3=\bar{\gamma}  ^n u_n.
\end{align}
\end{subequations}
Hence, the \lq prolongation\rq \;of the spanning vectors of the Lie algebra of \eqref{un} are as follows:
\begin{subequations}\label{3sym}
\begin{align}
X_1=  &(-1)^n u_n \frac{\partial\; }{\partial u_{n}}+(-1)^{n+1} u_{n+1}  \frac{\partial \quad }{\partial u_{n+1}}+(-1)^{n+2} u_{n +2}  \frac{\partial \quad }{\partial u_{n+2}}+ (-1)^{n+3} u_{n+3}  \frac{\partial\quad }{\partial u_{n+3}},\\
 &\nonumber\\
X_2= &\gamma^{n} u_{n} \frac{\partial\; }{\partial u_{n}}+\gamma^{n+1} u_{n+1} \frac{\partial \quad }{\partial u_{n+1}}+\gamma^{n+2} u_{n +2} \frac{\partial \quad }{\partial u_{n+2}}+\gamma^{n+3} u_{n+3}\frac{\partial\quad }{\partial u_{n+3}},\\&\nonumber\\
X_3=&\bar{\gamma}^{n} u_{n} \frac{\partial\; }{\partial u_{n}}+\bar{\gamma}^{n+1} u_{n+1} \frac{\partial \quad }{\partial u_{n+1}}+\bar{\gamma}^{n+2} u_{n +2} \frac{\partial \quad }{\partial u_{n+2}}+\bar{\gamma}^{n+3} u_{n+3} \frac{\partial\quad }{\partial u_{n+3}}.
\end{align}
\end{subequations}
Using $\xi _2$, the canonical coordinate \cite{JV} that linearizes \eqref{un} is
\begin{align}\label{cano}
S_n =\int\frac{du_n}{\xi _2(n,u_n)}=\frac{1}{\gamma ^n}\ln|u_n|.
\end{align}
Now, inspired by the form of the \lq final constraint\rq \eqref{rel}, we define the invariant function
\begin{align}\label{cano2}
\tilde{V}_n = \gamma^ nS_n+ \gamma ^{n+3}S_{n+3}
\end{align}
since
\begin{align}\label{inv}
X_1 \tilde{V}_n =X_2 \tilde{V}_n =X_3 \tilde{V}_n =0.
\end{align}
 For the sake of simplicity, we will be making use of the function $V_n$ defined as:
\begin{align}\label{vnexpo}
|{V}_n| =\exp\{  -\tilde{V}_n \},
\end{align}
that is to say,
$V_n =\pm {1}/({u_nu_{n+3}})$
but we shall utilize the plus sign:
\begin{align}\label{dam}
V_n =\frac{1}{u_nu_{n+3}}.
\end{align}
Substituting \eqref{vnexpo} into equation \eqref{un}, we obtain a first-order linear difference equation:
\begin{equation}\label{vn1}
V_{n+1}={A_nV_n}+ B_n
\end{equation}
whose solution solution in closed form is
\begin{align}\label{solvn}
V_{n}\quad=&V_0 \left(\prod\limits_{k_1=0}^{n-1}A_{k_1}\right) +\sum\limits_{l=0}^{n-1} \left(  B_{l}\prod\limits_{k_2=l+1}^{n-1}A_{k_2}\right).
\end{align}
From \eqref{cano}, \eqref{cano2} and \eqref{vnexpo}, we get
\begin{align}\label{solunk}
|u_n|=&  \exp\bigg((-1)^{n}c_1 + \bar{\gamma}^{n}c_2 + \gamma^{n}c_{3} + (-1)^{n}\sum_{k_1 = 0}^{n - 1}\frac{1}{3}(-1)^{-k_1}\ln \vert V_{k_1}\vert + \bar{\gamma}^{n}\sum_{k_2 = 0}^{n - 1} \frac{1}{3}\gamma^{k_2}\ln |V_{k_2}|\nonumber\\
&+ \gamma^{n}\sum_{k_3 = 1}^{n - 1} \frac{1}{3}\bar{\gamma^{k_3}}\ln |V_{k_3}|\bigg)           \nonumber    \\
=&   \exp\left\{ H_n+\frac{1}{3} \sum_{k = 0}^{n - 1} \left[(-1)^{n - k} + 2\text{Re}(\mathbb{H}(n,k)) \right]\ln |V_{k}|\right\}
\end{align}
where $\text{Re}(z)$ is the real part of $z$,
\begin{subequations}
\begin{align}
&H _n = (-1)^{n}c_1 + \bar{\gamma}^{n}c_2 + \gamma^{n}c_{3},\\
&\mathbb{H}(n,k) =\gamma ^n \bar{\gamma}^k.
 \end{align}
 \end{subequations}
 Recall that $\gamma =\exp{i\pi/3}$. It is easy to verify that
 \begin{subequations}\label{prop}
\begin{align}\label{Gammap}
 &\exp (H _{6n+j})=\exp (H_j) =|u_j|, \qquad \qquad \quad 0\leq j\leq 2,\\
 &\exp (H _{6n+j})=\exp (H_j) =|V_{j-3}||u_{j}|, \qquad \qquad \quad 3\leq j\leq 5,
 \end{align}
\begin{align}\label{gammap}
& \mathbb{H}(6n+j,k)= \mathbb{H}(j,6n+k)=\mathbb{H}(j,k),\\
&\mathbb{H}(n+3,k)= \mathbb{H}(n,k+3)=-\mathbb{H}(n,k).
\end{align}
\end{subequations}
Using \eqref{vnexpo}, together with properties \eqref{prop}, in \eqref{solunk}, we get
\begin{align}\label{solun}
|u_{6n+j}|=& \exp\left\{ H_j+\frac{1}{3} \sum_{k = 0}^{6n+j - 1} \left[(-1)^{6n+j - k} + 2\text{Re}(\mathbb{H}(j,k)) \right]\ln |V_{k}|\right\}\\
= & |u_j|\prod_{s=0}^{n-1}\left|\frac{V_{6n+j}}{V_{6n+j+3}}\right|.
 \end{align}
\par \noindent
Using \eqref{dam}, it can be shown that there is not need of the absolute values, thus
\begin{align*}
u_{6n + j} =  & u_j\prod_{s=0}^{n-1}\frac{V_{6n+j}}{V_{6n+j+3}} \\
           = & u_j \prod_{s=0}^{n-1}\frac{V_0 \left(\prod\limits_{k_1=0}^{6n + j -1}A_{k_1}\right) +\sum\limits_{l=0}^{6n + j-1} \left(  B_{l}\prod\limits_{k_2=l+1}^{6n + j -1}A_{k_2}\right)}{V_0 \left(\prod\limits_{k_1=0}^{6n + j + 2}A_{k_1}\right) +\sum\limits_{l=0}^{6n + j + 2} \left(  B_{l}\prod\limits_{k_2=l+1}^{6n + j + 2}A_{k_2}\right)}\\
           = & u_j \prod_{s=0}^{n-1}\frac{\left(\prod\limits_{k_1=0}^{6n + j -1}A_{k_1}\right) + \frac{1}{V_0}\sum\limits_{l=0}^{6n + j-1} \left(  B_{l}\prod\limits_{k_2=l+1}^{6n + j -1}A_{k_2}\right)}{\left(\prod\limits_{k_1=0}^{6n + j + 2}A_{k_1}\right) + \frac{1}{V_0}\sum\limits_{l=0}^{6n + j + 2} \left(  B_{l}\prod\limits_{k_2=l+1}^{6n + j + 2}A_{k_2}\right)}\\
    = & u_j \prod_{s=0}^{n-1}\frac{\left(\prod\limits_{k_1=0}^{6n + j -1}A_{k_1}\right) + u_0u_3\sum\limits_{l=0}^{6n + j-1} \left(  B_{l}\prod\limits_{k_2=l+1}^{6n + j -1}A_{k_2}\right)}{\left(\prod\limits_{k_1=0}^{6n + j + 2}A_{k_1}\right) + u_0u_3\sum\limits_{l=0}^{6n + j + 2} \left(  B_{l}\prod\limits_{k_2=l+1}^{6n + j + 2}A_{k_2}\right)}
\end{align*}
so that, for $j = 0,1,2,3,4,5$, we have
\begin{equation}\label{sol0}
x_{6n + j - 3} = x_{j - 3} \prod_{s=0}^{n-1}\frac{\left(\prod\limits_{k_1=0}^{6n + j -1}a_{k_1}\right) + x_{-3}x_0\sum\limits_{l=0}^{6n + j-1} \left(  b_{l}\prod\limits_{k_2=l+1}^{6n + j -1}a_{k_2}\right)}{\left(\prod\limits_{k_1=0}^{6n + j + 2}a_{k_1}\right) + x_{-3}x_0\sum\limits_{l=0}^{6n + j + 2} \left(  b_{l}\prod\limits_{k_2=l+1}^{6n + j + 2}a_{k_2}\right)}.
\end{equation}
Hence, the solution $\{x_{n}\}_{n = 1}^{\infty}$expressed in terms of the initial values $x_{-3}, x_{-2}, x_{-1}, x_{0}$ is given by;
\begin{equation}\label{sol1}
x_{6n - 3} = x_{- 3} \prod_{s=0}^{n-1}\frac{\left(\prod\limits_{k_1=0}^{6n-1}a_{k_1}\right) + x_{-3}x_0\sum\limits_{l=0}^{6n -1} \left(  b_{l}\prod\limits_{k_2=l+1}^{6n-1}a_{k_2}\right)}{\left(\prod\limits_{k_1=0}^{6n + 2}a_{k_1}\right) + x_{-3}x_0\sum\limits_{l=0}^{6n + 2} \left(  b_{l}\prod\limits_{k_2=l+1}^{6n + 2}a_{k_2}\right)},
\end{equation}
\begin{equation}\label{sol2}
x_{6n - 2} = x_{-2} \prod_{s=0}^{n-1}\frac{\left(\prod\limits_{k_1=0}^{6n}a_{k_1}\right) + x_{-3}x_0\sum\limits_{l=0}^{6n} \left(  b_{l}\prod\limits_{k_2=l+1}^{6n}a_{k_2}\right)}{\left(\prod\limits_{k_1=0}^{6n + 3}a_{k_1}\right) + x_{-3}x_0\sum\limits_{l=0}^{6n + 3} \left(  b_{l}\prod\limits_{k_2=l+1}^{6n + 3}a_{k_2}\right)},
\end{equation}
\begin{equation}\label{sol3}
x_{6n - 1} = x_{-1} \prod_{s=0}^{n-1}\frac{\left(\prod\limits_{k_1=0}^{6n + 1}a_{k_1}\right) + x_{-3}x_0\sum\limits_{l=0}^{6n + 1} \left(  b_{l}\prod\limits_{k_2=l+1}^{6n + 1}a_{k_2}\right)}{\left(\prod\limits_{k_1=0}^{6n + 4}a_{k_1}\right) + x_{-3}x_0\sum\limits_{l=0}^{6n + 4} \left(  b_{l}\prod\limits_{k_2=l+1}^{6n + 4}a_{k_2}\right)},
\end{equation}
\begin{equation}\label{sol4}
x_{6n} = x_{0} \prod_{s=0}^{n-1}\frac{\left(\prod\limits_{k_1=0}^{6n + 2}a_{k_1}\right) + x_{-3}x_0\sum\limits_{l=0}^{6n + 2} \left(  b_{l}\prod\limits_{k_2=l+1}^{6n + 2}a_{k_2}\right)}{\left(\prod\limits_{k_1=0}^{6n + 5}a_{k_1}\right) + x_{-3}x_0\sum\limits_{l=0}^{6n + 5} \left(  b_{l}\prod\limits_{k_2=l+1}^{6n + 5}a_{k_2}\right)},
\end{equation}
\begin{equation}\label{sol5}
x_{6n + 1} = \frac{x_{-3}x_0}{x_{-2}(a_0 + b_0x_{-3}x_0)} \prod_{s=0}^{n-1}\frac{\left(\prod\limits_{k_1=0}^{6n + 3}a_{k_1}\right) + x_{-3}x_0\sum\limits_{l=0}^{6n + 3} \left(  b_{l}\prod\limits_{k_2=l+1}^{6n + 3}a_{k_2}\right)}{\left(\prod\limits_{k_1=0}^{6n + 6}a_{k_1}\right) + x_{-3}x_0\sum\limits_{l=0}^{6n + 6} \left(  b_{l}\prod\limits_{k_2=l+1}^{6n + 6}a_{k_2}\right)},
\end{equation}
\begin{equation}\label{sol6}
x_{6n + 2} = \frac{x_{-3}x_0}{x_{-1}(a_0a_1 + (a_1b_0 + b_1)x_{-3}x_0)} \prod_{s=0}^{n-1}\frac{\left(\prod\limits_{k_1=0}^{6n + 4}a_{k_1}\right) + x_{-3}x_0\sum\limits_{l=0}^{6n + 4} \left(  b_{l}\prod\limits_{k_2=l+1}^{6n + 4}a_{k_2}\right)}{\left(\prod\limits_{k_1=0}^{6n + 7}a_{k_1}\right) + x_{-3}x_0\sum\limits_{l=0}^{6n + 7} \left(  b_{l}\prod\limits_{k_2=l+1}^{6n + 7}a_{k_2}\right)}.
\end{equation}
In the following section we look at a special case when the sequences $a_n$ and $b_n$ are 1-periodic.
\subsection{The case $a_n$ and $b_n$ are 1-periodic}
In this case, set $a_n = a$ and $b_n = b$ where $a, b\in \mathbb{R}$. Then by direct substitution in the solution above, we have\\
\begin{minipage}{.4\linewidth}
\begin{equation*}\label{soll1}
x_{6n - 3} = x_{- 3} \prod_{s=0}^{n-1}\frac{ a^{6n} + bx_{-3}x_0\sum\limits_{l=0}^{6n -1}a^{l}}{a^{6n + 3} + bx_{-3}x_0\sum\limits_{l=0}^{6n + 2}a^{l}},\;
\end{equation*}
\end{minipage}
\begin{minipage}{.4\linewidth}
\begin{equation*}\label{soll2}
\quad x_{6n - 2} = x_{-2} \prod_{s=0}^{n-1}\frac{ a^{6n + 1} + bx_{-3}x_0\sum\limits_{l=0}^{6n}a^{l}}{a^{6n + 4} + bx_{-3}x_0\sum\limits_{l=0}^{6n + 3}a^{l}},
\end{equation*}
\end{minipage}\\
\begin{minipage}{.4\linewidth}
\begin{equation*}\label{soll3}
x_{6n - 1} = x_{-1} \prod_{s=0}^{n-1}\frac{ a^{6n + 2} + bx_{-3}x_0\sum\limits_{l=0}^{6n + 1}a^{l}}{a^{6n + 5} + bx_{-3}x_0\sum\limits_{l=0}^{6n + 4}a^{l}},
\end{equation*}
\end{minipage}
\begin{minipage}{.4\linewidth}
\begin{equation*}\label{soll4}
 \quad x_{6n} = x_{0} \prod_{s=0}^{n-1}\frac{ a^{6n + 3} + bx_{-3}x_0\sum\limits_{l=0}^{6n + 2}a^{l}}{a^{6n + 6} + bx_{-3}x_0\sum\limits_{l=0}^{6n + 5}a^{l}},
\end{equation*}
\end{minipage}\\
\begin{minipage}{.4\linewidth}
\begin{equation*}\label{soll5}
x_{6n + 1} = c\prod_{s=0}^{n-1}\frac{ a^{6n + 4} + bx_{-3}x_0\sum\limits_{l=0}^{6n + 3}a^{l}}{a^{6n + 7} + bx_{-3}x_0\sum\limits_{l=0}^{6n + 6}a^{l}},
\end{equation*}
\end{minipage}
\begin{minipage}{.4\linewidth}
\begin{equation*}\label{soll6}
x_{6n + 2} = d\prod_{s=0}^{n-1}\frac{ a^{6n + 5} + bx_{-3}x_0\sum\limits_{l=0}^{6n + 4}a^{l}}{a^{6n + 8} + bx_{-3}x_0\sum\limits_{l=0}^{6n + 7}a^{l}},
\end{equation*}
\end{minipage}\\\\
where $c =  \frac{x_{-3}x_0}{x_{-2}(a + bx_{-3}x_0)}$ and $d = \frac{x_{-3}x_0}{x_{-1}(a^{2} + (ab + b)x_{-3}x_0)}$.
\subsubsection{The case $a = 1$}
\begin{minipage}{.4\linewidth}
\begin{equation*}\label{so1}
x_{6n - 3} = x_{- 3} \prod_{s=0}^{n-1}\frac{1 + (6n)bx_{-3}x_0}{1 + (6n + 3)bx_{-3}x_0},
\end{equation*}
\end{minipage}
\begin{minipage}{.4\linewidth}
\begin{equation*}\label{so2}
x_{6n - 2} = x_{-2} \prod_{s=0}^{n-1}\frac{1 + (6n + 1)bx_{-3}x_0}{1 + (6n + 4)bx_{-3}x_0},
\end{equation*}
\end{minipage}\\
\begin{minipage}{.4\linewidth}
\begin{equation*}\label{so3}
x_{6n - 1} = x_{-1} \prod_{s=0}^{n-1}\frac{1 + (6n + 2)bx_{-3}x_0}{1 + (6n + 5)bx_{-3}x_0},
\end{equation*}
\end{minipage}
\begin{minipage}{.4\linewidth}
\begin{equation*}\label{so4}
 x_{6n} = x_{0} \prod_{s=0}^{n-1}\frac{1 + (6n + 3)bx_{-3}x_0}{1 + (6n + 6)bx_{-3}x_0},
\end{equation*}
\end{minipage}\\
\begin{minipage}{.4\linewidth}
\begin{equation*}\label{so5}
x_{6n + 1} = c\prod_{s=0}^{n-1}\frac{1 + (6n + 4)bx_{-3}x_0}{1 + (6n + 7)bx_{-3}x_0},
\end{equation*}
\end{minipage}
\begin{minipage}{.4\linewidth}
\begin{equation*}\label{so6}
x_{6n + 2} = d\prod_{s=0}^{n-1}\frac{ 1 + (6n + 5)bx_{-3}x_0}{1 + (6n + 8)bx_{-3}x_0}.
\end{equation*}
\end{minipage} \\\\
where $c = \frac{x_{-3}x_0}{x_{-2}(1 + bx_{-3}x_0)}$ and $d = \frac{x_{-3}x_0}{x_{-1}(1 + 2bx_{-3}x_0)}$.
\subsubsection{The case $a \neq 1$}
For this case, we obtain\\
\begin{minipage}{.4\linewidth}
\begin{equation*}\label{soll1}
x_{6n - 3} = x_{- 3} \prod_{s=0}^{n-1}\frac{ a^{6n} + bx_{-3}x_0\left(\frac{1-a^{6n}}{1-a} \right)}{a^{6n + 3} + bx_{-3}x_0\left(\frac{1-a^{6n+3}}{1-a} \right)},\;
\end{equation*}
\end{minipage}
\begin{minipage}{.4\linewidth}
\begin{equation*}\label{soll2}
\qquad x_{6n - 2} = x_{-2} \prod_{s=0}^{n-1}\frac{ a^{6n + 1} + bx_{-3}x_0\left(\frac{1-a^{6n+1}}{1-a} \right)}{a^{6n + 4} + bx_{-3}x_0\left(\frac{1-a^{6n+4}}{1-a} \right)},
\end{equation*}
\end{minipage}\\
\begin{minipage}{.4\linewidth}
\begin{equation*}\label{soll3}
x_{6n - 1} = x_{-1} \prod_{s=0}^{n-1}\frac{ a^{6n + 2} + bx_{-3}x_0\left(\frac{1-a^{6n+2}}{1-a} \right)}{a^{6n + 5} + bx_{-3}x_0\left(\frac{1-a^{6n+5}}{1-a} \right)},
\end{equation*}
\end{minipage}
\begin{minipage}{.4\linewidth}
\begin{equation*}\label{soll4}
 \qquad \quad x_{6n} = x_{0} \prod_{s=0}^{n-1}\frac{ a^{6n + 3} + bx_{-3}x_0\left(\frac{1-a^{6n+3}}{1-a} \right)}{a^{6n + 6} + bx_{-3}x_0\left(\frac{1-a^{6n+6}}{1-a} \right)},
\end{equation*}
\end{minipage}\\
\begin{minipage}{.4\linewidth}
\begin{equation*}\label{soll5}
x_{6n + 1} = c\prod_{s=0}^{n-1}\frac{ a^{6n + 4} + bx_{-3}x_0\left(\frac{1-a^{6n+4}}{1-a} \right)}{a^{6n + 7} + bx_{-3}x_0\left(\frac{1-a^{6n+7}}{1-a} \right)},
\end{equation*}
\end{minipage}
\begin{minipage}{.4\linewidth}
\begin{equation*}\label{soll6}
\qquad x_{6n + 2} = d\prod_{s=0}^{n-1}\frac{ a^{6n + 5} + bx_{-3}x_0\left(\frac{1-a^{6n+5}}{1-a} \right)}{a^{6n + 8} + bx_{-3}x_0\left(\frac{1-a^{6n+8}}{1-a} \right)},
\end{equation*}
\end{minipage}\\\\
where $c =  \frac{x_{-3}x_0}{x_{-2}(a + bx_{-3}x_0)}$ and $d = \frac{x_{-3}x_0}{x_{-1}(a^{2} + (ab + b)x_{-3}x_0)}$.\par
\vspace{1cm}
\noindent \underline{\textit{The case $a = - 1$}:}\\
For this case, we obtain the following solution:
$$x_{6n + j - 3}
= \begin{cases}
x_{j - 3}(-1 + bx_{-3}x_0)^{n}, & \text{if}\,\,j\,\,\text{is odd}; \\
x_{j - 3}(-1 + bx_{-3}x_0)^{-n}, & \text{if}\,\,j\,\,\text{is even};
\end{cases}
$$
where $j \in \{ 0, 1,2,3,4,5\}$, $x_{1} =  \frac{x_{-3}x_0}{x_{-2}(-1 + bx_{-3}x_0)}$ and $x_2 = \frac{x_{-3}x_0}{x_{-1}}$.
\section{Conclusion}
In this paper, we obtained three non-trivial symmetry generators and formulas for the solutions of difference equations \eqref{xn}. Our approach involved Lie symmetry analysis and solving certain recurrences relations.

\end{document}